\documentclass[11pt]{amsart} 

\usepackage{euscript}
\usepackage{amsmath}
\usepackage{amsthm}
\usepackage{epsfig}
\usepackage{amssymb}\usepackage{amscd}
\usepackage{epic}
\usepackage{eufrak}
\usepackage{tikz,underscore}
 
\numberwithin{equation}{section}
\theoremstyle{plain}
\newtheorem{theorem}[section]{Theorem}
\newtheorem{lemma}[section]{Lemma}

\theoremstyle{definition}
\newtheorem{example}[section]{Example}
\newtheorem{remark}[section]{Remark}
\newtheorem{definition}[section]{Definition}

\title{Links of plane curve singularities are $L$-space links}

\author{Eugene Gorsky}
\address{Department of Mathematics, Columbia University, 2990 Broadway, New York NY, USA}
\email{egorsky@math.columbia.edu}

\author{Andr\'as N\'emethi}
\address{A. R\'enyi Institute of Mathematics, 1053 Budapest, Re\'altanoda u. 13-15, Hungary.}
\email{nemethi.andras@renyi.mta.hu}

\thanks{\noindent The first author is partially supported by the grants
RFBR–13-01-00755, NSh-4850.2012.1.
The second author is partially supported by OTKA Grant 100796.}
 
\begin{document}

\begin{abstract}
 We prove that a sufficiently large surgery on any algebraic link is an $L$-space.
 For torus links we give a complete classification of integer surgery coefficients 
 providing $L$-spaces. 
\end{abstract}

\maketitle

\begin{definition}
A 3-manifold $Y$ is called an $L$-space, if it is a rational homology sphere and  its Heegaard-Floer homology has minimal possible rank:
$\mathrm{rk}\ \widehat{HF}(Y)=|H_1(Y,\mathbb{Z})|.$
\end{definition}

We refer the reader to \cite{os2,os3,os4} for definitions of Heegaard-Floer homology, and to \cite{hedden,OS}  for the detailed discussion of the properties of $L$-spaces. Let $K\subset S^3$ be the embedded  link of a  complex plane curve singularity with $r$ components:
$K=K_1\cup\cdots \cup K_r.$ The following theorem is the main result of this note.

\begin{theorem}
\label{th1}
Every algebraic link $K\subset S^3$ is an $L$--space link.
This means that an integral surgery of $S^3 $ along the link components $K_i$
with all coefficients sufficiently large  provides an $L$--space.
\end{theorem}

For $r=1$ the result was proved by Hedden \cite[Theorem 1.10]{hedden}.
Our proof is of different nature, is extremely short,
and provides the argument uniformly for any $r$. It is based on some facts from theory of normal surface singularities and plumbed 3-manifolds.
 For plumbing calculus (which modifies the possible graph representations of the same 3-manifold)
 see \cite{neumann}.
We will need the following facts:

(a) A connected negative definite  plumbing graph is called a ``smooth graph'' if by blowing down
consecutively
$(-1)$--vertices we can blow down the graph to the empty graph. Such a graph represents $S^3$.
When we resolve plane curve singularities and we blow up $({\mathbb C}^2,0)$ in several infinitely
near points we obtain such a graph.

(b) A surface singularity is rational if its geometric genus is zero. This property can be verified at the level
of its negative definite plumbing graph (such graphs  are called ``rational graphs''),
see \cite{Artin62,Artin66,Laufer72}.

Blowing up a rational graph we get a rational graph.
Any subgraph of a rational graph is rational (use e.g. Laufer's
criterion \cite{Laufer72}).  Since smooth graphs are
rational, subgraphs of smooth graphs are rational (they are called ``sandwiched'' graphs, see
\cite{S}).  The following remark describes two useful consequences of Laufer's criterion. 

\begin{remark}(cf. \cite[Rmk 2.3, Prop 2.4]{S})
\label{laufer}
For each vertex $v$ of the negative definite plumbing graph $\Gamma$ define $w(v)=-(E_v,E_v)$, where $E_v$ is the corresponding curve.
Let $\gamma(v)$ denote the valency of $v$ in $\Gamma$.

(a) If $\Gamma$ is rational, then $w(v)\ge \gamma(v)-1$ for all $v$. 

(b)  If $w(v)\ge \gamma(v)$ for all $v$, then the $\Gamma$ is rational (and, in fact, sandwiched).
Such graphs are called minimal rational.
\end{remark}

A key ingredient in the proof of Theorem \ref{th1} is the following result:

\begin{theorem}\cite[Thms 6.3, 8.3]{N} A  3--manifold
 plumbed from a rational graph is an $L$--space.
 \end{theorem}

%\vspace{2mm}

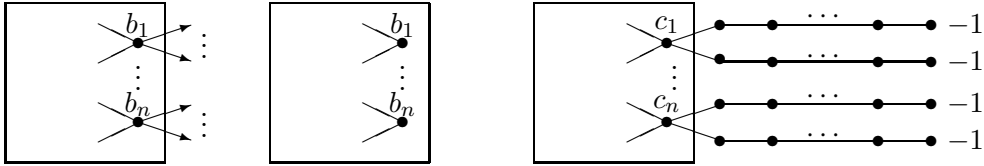
\begin{figure}[ht]

\begin{picture}(100,80)(150,0)
\put(10,10){\framebox(60,60)}
\put(60,25){\circle*{4}} \put(60,55){\circle*{4}}
\put(60,32){\makebox(0,0){$b_n$}}
\put(60,62){\makebox(0,0){$b_1$}}
\put(60,25){\line(-2,1){15}}  \put(60,25){\line(-2,-1){15}}
\put(60,55){\line(-2,1){15}}  \put(60,55){\line(-2,-1){15}}
\put(60,45){\makebox(0,0){$\vdots$}}
\put(60,25){\vector(3,1){20}}
\put(60,55){\vector(3,1){20}}
\put(60,25){\vector(3,-1){20}}
\put(60,55){\vector(3,-1){20}}
\put(85,27){\makebox(0,0){$\vdots$}}
\put(85,57){\makebox(0,0){$\vdots$}}

\put(110,10){\framebox(60,60)}
\put(160,25){\circle*{4}} \put(160,55){\circle*{4}}
\put(160,32){\makebox(0,0){$b_n$}}
\put(160,62){\makebox(0,0){$b_1$}}
\put(160,25){\line(-2,1){15}}  \put(160,25){\line(-2,-1){15}}
\put(160,55){\line(-2,1){15}}  \put(160,55){\line(-2,-1){15}}
\put(160,45){\makebox(0,0){$\vdots$}}

\put(210,10){\framebox(60,60)}
\put(260,25){\circle*{4}} \put(260,55){\circle*{4}}
\put(260,32){\makebox(0,0){$c_n$}}
\put(260,62){\makebox(0,0){$c_1$}}
\put(260,25){\line(-2,1){15}}  \put(260,25){\line(-2,-1){15}}
\put(260,55){\line(-2,1){15}}  \put(260,55){\line(-2,-1){15}}
\put(263,45){\makebox(0,0){$\vdots$}}
\put(260,25){\line(3,1){20}}
\put(260,25){\line(3,-1){20}}
\put(260,55){\line(3,1){20}}
\put(260,55){\line(3,-1){20}}
\put(280,32){\line(1,0){80}}
\put(280,18){\line(1,0){80}}
\put(280,62){\line(1,0){80}}
\put(280,48){\line(1,0){80}}

\put(300,32){\circle*{4}} \put(300,62){\circle*{4}}
\put(340,32){\circle*{4}} \put(340,62){\circle*{4}}
\put(360,32){\circle*{4}} \put(360,62){\circle*{4}}
\put(280,32){\circle*{4}} \put(280,62){\circle*{4}}
\put(280,18){\circle*{4}} \put(280,49){\circle*{4}}
\put(300,18){\circle*{4}} \put(300,48){\circle*{4}}
\put(340,18){\circle*{4}} \put(340,48){\circle*{4}}
\put(360,18){\circle*{4}} \put(360,48){\circle*{4}}
\put(320,35){\makebox(0,0){$\cdots$}}\put(320,65){\makebox(0,0){$\cdots$}}
\put(320,20){\makebox(0,0){$\cdots$}}\put(320,50){\makebox(0,0){$\cdots$}}
 
\put(373,18){\makebox(0,0){$-1$}}
\put(373,32){\makebox(0,0){$-1$}}
\put(373,48){\makebox(0,0){$-1$}}
\put(373,62){\makebox(0,0){$-1$}}

\end{picture}
\caption{Graphs $\Gamma_K,\Gamma_0$ and $\widetilde{\Gamma_0}$.}
\label{fig1}
\end{figure}

Let $C=C_1\cup \ldots \cup C_r\subset \mathbb{C}^2$ be the plane curve singularity, corresponding to $K$,
so that the link of the irreducible component $C_i$  is the knot $K_i\subset S^3$.
The possible minimal embedded resolution (plumbing) graphs of complex plane curve singularities are
well--known.  
We will represent such a graph in a schematic way by
the  graph $\Gamma_K$ shown in Figure \ref{fig1}, emphasizing only those
exceptional curves, say $E_1,\ldots, E_n$ which intersect one or more strict transform components.
The strict transform components are encoded by arrowheads. The number of arrowheads supported by
$E_i$ is $a_i$, and the self-intersection number of $E_i$ is $b_i$ ($1\leq i\leq n$). Definitely, some of  $b_i$'s are equal to $-1$.

If we delete the arrowheads, we get a smooth graph  $\Gamma_0$ (the second graph in Figure \ref{fig1}).
If we blow up $\Gamma_0$ few times (starting with $E_i$'s and continuing with the newly created
$(-1)$ vertices),
we obtain another graph $\widetilde{\Gamma_0}$,
which is again smooth.  This is the third  graph in Figure \ref{fig1}.
The unmarked vertices are $(-2$)--vertices and  $c_i=b_i-a_i$.
The length of the newly created legs can be different:
the number of $(-2)$-vertices in the $(i,j)$th leg
is $(k_{ij}-1)$, where $k_{ij}\geq 1$ for all  $1\leq i\leq n$, $1\leq j\leq a_i$.

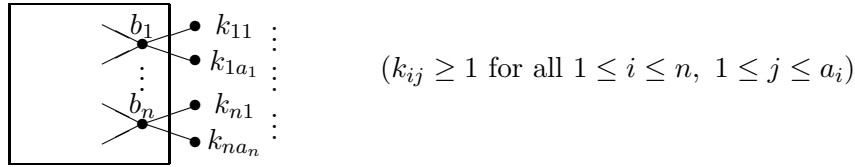
\begin{figure}[ht]
\begin{picture}(100,70)(100,0)
\put(10,10){\framebox(60,60)}
\put(60,25){\circle*{4}} \put(60,55){\circle*{4}}
\put(60,32){\makebox(0,0){$b_n$}}
\put(60,62){\makebox(0,0){$b_1$}}
\put(60,25){\line(3,1){20}}
\put(60,25){\line(3,-1){20}}
\put(60,55){\line(3,1){20}}
\put(60,55){\line(3,-1){20}}
\put(80,32){\circle*{4}} \put(80,62){\circle*{4}}
\put(80,18){\circle*{4}} \put(80,49){\circle*{4}}

\put(60,25){\line(-2,1){15}}  \put(60,25){\line(-2,-1){15}}
\put(60,55){\line(-2,1){15}}  \put(60,55){\line(-2,-1){15}}
\put(60,45){\makebox(0,0){$\vdots$}}
 
 \put(95,32){\makebox(0,0){$k_{n1}$}}
\put(95,62){\makebox(0,0){$k_{11}$}}
 \put(95,18){\makebox(0,0){$k_{na_n}$}}
\put(95,48){\makebox(0,0){$k_{1a_1}$}}
\put(110,60){\makebox(0,0){$\vdots$}}\put(110,45){\makebox(0,0){$\vdots$}}
\put(110,27){\makebox(0,0){$\vdots$}}
\put(150,45){\makebox(0,0)[l]{($k_{ij}\geq 1$ for all $1\leq i\leq n, \  1\leq j\leq a_i$)}}
\end{picture}
\caption{Graph $\Gamma$}
\label{fig2}
\end{figure}

\begin{lemma}
\label{lem1}
For $k_{ij}\ge 1$  the 3-manifold represented by the graph $\Gamma$ shown in Figure \ref{fig2} is an $L$-space.
\end{lemma}

\begin{proof}
By plumbing calculus (blowing up several times the most right edges)
$\Gamma$ is equivalent to the graph $\Gamma'$
obtained from  the smooth graph $\widetilde{\Gamma_0}$
by deleting the $(-1)$ vertices from the end of the new legs.
Hence $\Gamma'$  is a subgraph of a smooth graph, thus it is rational.
In particular, it represents an $L$--space.
\end{proof}
\noindent
{\bf Proof of Theorem \ref{th1}:} Suppose that the plane curve singularity, corresponding to $K$, is given by the equation
$\{\prod_{i,j} f_{ij}(x,y)=0\}$.
Consider a positive surgery on $S^3$ along the link $K$, where the surgery coefficient for the component $K_{ij}$ equals $d_{ij}$.
It is known that $S^3_{\{d_{ij}\}_{ij}}(K)$ is a plumbed 3-manifold with the plumbing graph
 $\Gamma$, such that
the parameters $k_{ij}$ are defined by the equation $k_{ij}=d_{ij}-m_i$, where $m_i$
 are the multiplicities of the pullback of $f_{ij}$ on the divisor $E_i$.
For $d_{ij}>m_i$ the surgery space is an $L$-space by Lemma \ref{lem1}. $\square$

\vspace{2mm}

\begin{remark} \ (a) Note that the above $L$-spaces are very special: they are given by
connected negative definite sandwiched graphs (e.g, they are irreducible),
so they definitely have even  more additional properties.

(b) The above bound $(k_{ij}\geq 1)$ is not optimal, usually one can find a
collection of smaller numbers $B_{ij}$, $B_{ij}$ smaller than 1,
such that all the surgeries with  coefficients $k_{ij}\geq B_{ij}$ 
provide  $L$--spaces. % (see the bounds in Hedden \cite{hedden}).
But for the coefficients $B_{ij}=1$ the proof is extremely transparent
(and for smaller coefficients singularity theory is harder to apply).

(c) If $a_i=1$ for all $i$ then one can take  $k_{i1}\geq 0$, and the surgery manifold is still
  an $L$--space. For example, if we take all of them zero
  then $\Gamma$ is equivalent to a graph which is obtained from $\Gamma$ by deleting the 0
  vertices and the
  supporting $E_i$ vertices.
This graph is  not connected, but each component
 is a subgraph of a smooth graph. Hence the corresponding 3-manifold is a
connected sum of $L$--spaces, which is again an $L$-space.
The proof of the general case $k_{i1}\geq 0$ (with $a_i=1$) is a combination of this argument with
the proof of Theorem 2.

(d) It is known that all algebraic knots (with one component) can be presented as iterated cables of the trivial knot. Hedden proved in \cite{hedden} that the $pq$ surgery of $S^3$ along an algebraic knot is an $L$-space, where $p$ and $q$ are the parameters of the last cabling. One can check that $pq=m_1$ in this case.  
See also \cite{hom} for a complete description of $L$-space surgeries of cable knots.
\end{remark}

It turns out that the set of all surgeries on a link providing an $L$-space has an interesting structure.
This is new phenomenon compared with the irreducible case, where (e.g. by \cite{hom}) 
a $d$-surgery on a knot $K$ is an $L$-space iff $d\ge 2g(K)-1$ (hence $d$ runs in a half-line).
The following theorem provides a  description of $L$-space surgery coefficients for torus links.

%\begin{example}
%The coefficient $k_{ij}=0$ has an interesting peculiarity when $a_i>1$.
%We exemplify it on the graph of the curve singularity $(\{x^6+y^9=0\},0)$:

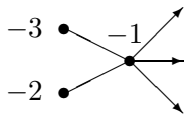
\begin{figure}[ht]
\begin{picture}(100,40)(0,20)
\put(20,28){\circle*{4}}
\put(20,52){\circle*{4}}\put(45,40){\circle*{4}}
\put(45,40){\line(-2,-1){25}}\put(45,40){\line(-2,1){25}}
\put(45,40){\vector(1,1){20}}\put(45,40){\vector(1,0){20}}
\put(45,40){\vector(1,-1){20}}
\put(5,28){\makebox(0,0){$-2$}}\put(5,52){\makebox(0,0){$-3$}}
\put(43,50){\makebox(0,0){$-1$}}
\end{picture}
\caption{Plumbing graph for  the singularity $x^6+y^9=0$.}
\label{fig4}
\end{figure}

\begin{theorem}
\label{th2}
Let $p,q>1$ be two coprime integers, $r\ge 1$ and $(d_1,\ldots,d_r)\in \mathbb{Z}^r$.
Assume that $d_i\neq pq\pm 1$ for all $i$. Then $(d_1,\ldots,d_r)$--surgery on $(pr,qr)$--torus link
is an $L$-space iff the surgery space is a rational homology sphere and one of the following conditions hold:
\begin{itemize}
\item[(a)] $d_i=pq$ for some $i$. 
\item[(b)] $d_i>pq+1$ for all $i$.
\item[(c)] $d_i<pq-1$ for all $i$ and $\max(d_i)\ge pq-p-q$.
\end{itemize}
\end{theorem}

\begin{remark}
It is easy to see that $d_i=d_j=pq$ for $i\neq j$ yields infinite $H_1$ for the surgery space, so in Theorem \ref{th2}(a) 
$d_i=pq$ for exactly one $i$.
\end{remark}

\begin{remark}
The surgeries with $d_i=pq\pm 1$ can be also analyzed by (rather long) case analysis. In particular:

(a) The $(pq\pm 1,d)$--surgery on $(2p,2q)$--torus link is an $L$-space for all $d$.

(b) The $(pq+1,pq-1,d_3,\ldots,d_r)$--surgery on $(pr,qr)$--torus link coincides with $(d_3,\ldots,d_r)$ surgery
on $(p(r-2),q(r-2))$--torus link. In particular, $(pq+1,pq-1,d)$--surgery on a $(3p,3q)$--torus link is an $L$-space iff
$d\ge pq-p-q$.
\end{remark}

\begin{proof}
The plumbing graph of the surgery space can be described as above, with $k_i=d_i-pq$ (see Figure \ref{fig4} for the case of $(6,9)$ torus link with three components, see also \cite{torussurg} for more details). It is star-shaped, so the surgery space  $Y:=S^3_{d}(L)$ is Seifert fibered. 
By \cite{neumann2}, the intersection matrix for the Neumann normal form is negative definite  either for $Y$ or for $-Y$.
Any negative definite Seifert graph is almost rational in the sense of \cite{N}, so it represents an $L$-space if and only if it is rational.

If $d_i=pq$ then $k_i=0$ and the surgery space is a connected sum of lens spaces. %If $d_i=pq\pm 1$ for some $i$, one can blow down these vertices and check that the resulting graph is rational. 
From now on we will assume that $k_i\notin \{-1,0,1\}$ for all $i$. 

If $k_i>1$ for all $i$, then we get an $L$-space by Lemma \ref{lem1}. Suppose that $k_i$ have different signs, for example, $k_1>0$ and $k_2<0$. To obtain the normal form of $Y$, one needs to blow up the vertices with positive $k_i$, and there are at most $r-1$ of them, so the self-intersection of the central vertex in the normal form is greater than or equal to $-1-(r-1)=-r$, while its valency is $r+2$. To obtain the normal form of $-Y$, one reverses all signs and blows up all positive vertices (at most $r+1$ of them), so the self-intersection of the central vertex is greater than or equal to $1-(r+1)=-r$. If $Y$ is an $L$-space, either $Y$ or $-Y$ should be a negative definite rational graph, what contradicts Remark \ref{laufer}(a).

Finally, suppose that $k_i<0$ for all $i$.  In this case we can use a theorem of Lisca and Stipsicz \cite{StLi} describing Seifert fibered $L$-spaces.
Suppose that $1\ge \alpha_1\ge \alpha_2\ge \ldots \ge \alpha_{r+2}\ge 0$ are the Seifert invariants of singular fibers, and the central vertex has self-intersection $(-1)$. Then the corresponding 3-manifold is an $L$-space if and only if there are no coprime integers $(l,m)$ such that
\begin{equation}
\label{seifert l space}
m\alpha_1<l<m(1-\alpha_2)\ \text{and}\ m\alpha_i<1\ \text{for}\ i\ge 3.
\end{equation}
One can check that in our situation $\alpha_1$ and $1-\alpha_2$ are two neighboring fractions (in the sense of Farey series) with denominators $p$ and $q$ (that is, $\frac{a}{p}$ and $\frac{b}{q}$ with $aq-bp=\pm 1$), and $\alpha_i=1/|k_{i-2}|$ for $i\ge 3$. 

It is well known (see e.g. \cite[Section 4.5]{concrete}) that the fraction $\frac{a+b}{p+q}$ has the least possible denominator among rational numbers between  $\frac{a}{p}$ and $\frac{b}{q}$. Therefore if $k_i\ge -p-q$ for some $i$ then \eqref{seifert l space} cannot be satisfied and $Y$ is an $L$-space;
if $k_i<-p-q$ for all $i$ then $(l,m)=(a+b,p+q)$ satisfies \eqref{seifert l space} and $Y$ is not an $L$-space.
\end{proof}

\begin{example}
For $r=1$ the conditions of Theorem \ref{th2} are equivalent to the inequality $d_1\ge pq-p-q=2g(T(p,q))-1$, which also follows from \cite{hom}.
\end{example}

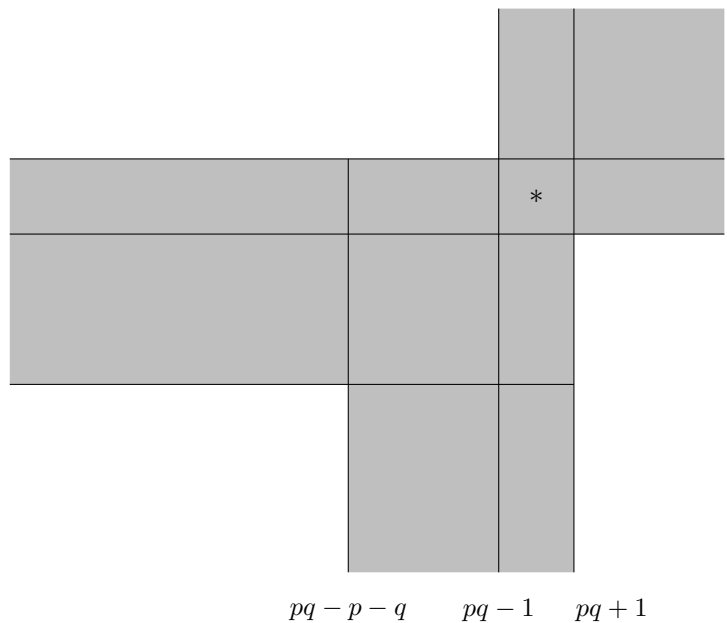
\begin{figure}
\begin{tikzpicture}
\fill [fill=lightgray] (-2.5,-5)--(0.5,-5)--(0.5,-0.5)--(2.5,-0.5)--(2.5,2.5)--(-0.5,2.5)--(-0.5,0.5)--(-7,0.5)--(-7,-2.5)--(-2.5,-2.5)--(-2.5,-5);
\draw (-2.5,-5)--(-2.5,0.5);
\draw (0.5,-5)--(0.5,2.5);
\draw (-0.5,-5)--(-0.5,2.5);
\draw (-7,-0.5)--(2.5,-0.5);
\draw (-7,0.5)--(2.5,0.5);
\draw (-7,-2.5)--(0.5,-2.5);
\draw (0,0) node {$*$};
\draw (-2.5,-5.5) node {\small $pq-p-q$};
\draw (-0.5,-5.5) node {\small $pq-1$};
\draw (1,-5.5) node {\small $pq+1$};
\end{tikzpicture}
\caption{Possible $L$-space surgery coefficients for $(2p,2q)$ torus link}
\label{L2}
\end{figure}

\begin{example}
Consider the case $r=2$.  
The set of all pairs $(d_1,d_2)$ providing an $L$-space is shown in grey in Figure \ref{L2}. The point $(pq,pq)$ (marked by a star) provides a surgery with infinite first homology and should be excluded. 
\end{example}

Theorem \ref{th1} can be reinterpreted as follows: start with a negative definite non-minimal
plumbing graph $\Gamma_0$
of $S^3$, put some arrows on the vertices, % (anywhere and as many as you wish),
and regard them as link components of $K$ in $S^3$. Then the surgery manifold with all
sufficiently large surgery coefficients is an $L$--space.

The proof works totally unmodified if we start with any non--minimal rational graph of a 3--manifold
$M$ (instead of $S^3$). If we wish to have natural (well--defined) surgery properties, it is convenient
to consider only the case when $M$ is an integral homology sphere. On the other hand,
this is rather restrictive for (negative definite) rational graphs, there are only two
possibilities: $M=S^3$ (case treated above),
and $M=\Sigma(2,3,5)$, the Poincar\'e 3--sphere oriented as the link of the surface singularity $\{x^2+y^3+z^5=0\}$.  Nevertheless, in this case, by identical proof we obtain:
%\end{remark}

\begin{theorem}
Consider any  negative definite non-minimal
plumbing graph of $\Sigma=\Sigma(2,3,5)$,
put some arrows on the vertices, % (anywhere and as many as you wish),
and regard them as link components of $K$ in $\Sigma$. Then the surgery manifold
$\Sigma(K)$ for all sufficiently large 
surgery coefficients is an $L$--space.
\end{theorem}

If we drop the integral homology restriction, we can start with any non--minimal rational graph
(now, this family is really large), but we have to consider only those surgeries
which can be realizes by the above construction (steps $\Gamma_0\mapsto \widetilde{\Gamma_0}\mapsto
\Gamma'$). For these surgeries the statement and the proof still holds.

\section*{Acknowledgments}

The authors would like to thank M. Hedden, J. Hom, R. Lipshitz, C. Manolescu, P. Ozsv\'ath and Z. Szab\'o for useful discussions. 
Special thanks to Jonathan Hanselman for sharing his program \cite{hanselman2} (based on \cite{hanselman}) which helped us
to formulate Theorem \ref{th2} and verify it in many examples.

\end{document}